\documentclass[11pt]{article}
\usepackage[T1]{fontenc}
\usepackage{amsmath,amssymb,amsthm,booktabs}
\usepackage[margin=1in]{geometry}
\usepackage[hidelinks]{hyperref}

\newtheorem{theorem}{Theorem}
\newtheorem{lemma}[theorem]{Lemma}
\newcommand{\R}{\mathbb{R}}
\newcommand{\Vone}{\mathrm{V1}}
\newcommand{\Vtwo}{\mathrm{V2}}

\title{Face-neighbourly polyhedra with eight or nine faces\\
and genus at most three}
\author{Ruslan Mizhaev\\[3pt]\small\href{mimodart@gmail.com}{mimodart@gmail.com}}
\date{}

\begin{document}
\maketitle

\begin{abstract}
A polyhedral surface is face-neighbourly if every pair of faces shares an edge.
We classify the combinatorial types of closed orientable face-neighbourly
polyhedral surfaces in \(\R^3\) with eight or nine faces and genus at most three,
assuming trivalent vertices, simple planar polygonal faces and distinct
supporting planes for adjacent faces. Several common edges between a pair of
faces are allowed. The enumeration gives \(42\,124\) candidate maps with eight
faces and \(203\) with nine. SAT certificates exclude all but two maps, both of
genus three with eight nonagonal faces. They are realised by Mizhaev's \(\Vtwo\)
and by the R\"ost--V\'igh polyhedron. With the face-neighbourliness condition
removed, the class of genus-three surfaces with eight nonagonal faces has one
additional realisable type, Mizhaev's \(\Vone\).
\end{abstract}

\section{Introduction}

Every two faces of the seven-faced Szilassi polyhedron share an edge
\cite{szilassi}. With eight faces, the same adjacency condition can be met if
some pairs are allowed to share several edges. Such intersections also occur
in the toroidal complexes studied by Gr\"unbaum and Szilassi \cite{gs}.

The 2020 preprint \cite{mizhaev2020} gives two genus-three polyhedra with eight
nonagonal faces, denoted here by their original labels \(\Vone\) and \(\Vtwo\).
Of these, \(\Vtwo\) is face-neighbourly. An integer-coordinate realisation of
\(\Vtwo\) appeared in \cite{mizhaev2026}. R\"ost and V\'igh subsequently
constructed a second face-neighbourly type with the same face sizes and genus
\cite{rostvigh}. Both types have eight pairs of faces that share two edges.
The graph of these pairs is an eight-cycle for \(\Vtwo\) and a disjoint union
of two four-cycles for the R\"ost--V\'igh example.

These constructions leave two questions. Do any other face-neighbourly types
with eight or nine faces occur through genus three? How many of the
eight-nonagon types occur if face-neighbourliness is dropped?

\section{Scope and results}

We consider finite, connected, closed, orientable polyhedral surfaces embedded
in \(\R^3\). Each face is a closed disk bounded by a simple planar polygon.
Face intersections consist of the vertices and edges prescribed by the map.
Every vertex has degree three and adjacent faces have distinct supporting
planes. We count the underlying cell decompositions up to isomorphism,
allowing reversal of orientation. A face-size profile is the multiset of the
numbers of sides. For example, \(7^4 8^4\) denotes four heptagons and four
octagons.

\begin{theorem}\label{thm:classification}
In this class the face-neighbourly surfaces with \(F\in\{8,9\}\) and genus
\(g\leq 3\) have exactly two combinatorial types. Both have \(F=8\), \(g=3\)
and profile \(9^8\). They are \(\Vtwo\) and the R\"ost--V\'igh type.
Among surfaces of genus three with profile \(9^8\), the realisable types are
\(\Vone\), \(\Vtwo\) and the R\"ost--V\'igh type.
\end{theorem}

\section{Enumeration of maps}

Write \(V\), \(E\) and \(F\) for the numbers of vertices, edges and faces.
Trivalence and Euler's formula give
\[
  E=3(F-2+2g),\qquad V=2(F-2+2g).
\]
Face-neighbourliness requires \(E\geq\binom F2\). Hence the only
face-neighbourly cases to enumerate are \((F,g)=(8,2),(8,3),(9,3)\).

At each vertex, the link has three edges and gives three distinct incident
faces. Any two of them are adjacent there. Their three supporting planes
meet transversely. Indeed, if they met in a line, all three incident edges
would lie on that line, and two would overlap near the vertex. Thus the
three faces determine the vertex uniquely. In particular the same face
triple cannot label two vertices.

Let \(m_{ij}\) be the number of common edges of faces \(i,j\). The endpoints
of these edges use distinct third faces, so
\(m_{ij}\leq c=\lfloor(F-2)/2\rfloor\). In the face-neighbourly cases we
reserve one edge for each pair. For a face of size \(p_i\), its residual
degree is \(d_i=p_i-(F-1)\). The profiles are generated from the side sum
\(\sum_i p_i=2E\) and the bounds
\[
 F-1\le p_i\le \min\{V,(F-1)c,F-1+E-\tbinom F2\}.
\]
Counting residual incidences gives an additional bound for every subset
\(S\) of \(k\) faces.
\[
 \sum_{i\in S} d_i \le (c-1)k(k-1)
 +\sum_{j\notin S}\min\{d_j,(c-1)k\}.
\]
The two terms bound the contributions from edges within \(S\) and from
edges leaving \(S\). For each \(k\) it suffices to test the \(k\) largest
residual degrees. They maximise the left side and minimise the right.
For \(F=8\), \(144\) profiles pass the side-count
bounds and \(47\) are excluded by the incidence test. For \(F=9\), the
corresponding counts are \(7\) and \(1\).

For each remaining profile, enumerate sets of unordered triples of face
labels, with label \(i\) occurring \(p_i\) times. For each pair \(i,j\) the
triples containing that pair must occur in even number. In the
face-neighbourly cases that number is positive. Enumerate all perfect
matchings on these triples for each pair. A matched pair supplies the
endpoints of an edge. Every map in the class supplies one such choice of
triples and matchings, since each occurrence of \(i,j\) at a vertex is
the endpoint of exactly one edge shared by the two faces.

Retain pairings for which the boundary of each face is a single simple
cycle and the resulting map is connected, closed and orientable. Two
faces with the same vertex set are discarded. In a realisation such faces
would be adjacent by trivalence, and their common non-collinear vertices
would force their supporting planes to coincide. To identify isomorphic
maps, canonically label the directed edge occurrences using the successor
around a face and reversal across an edge. Take the lexicographic minimum
over all starting occurrences and both global orientations. These two
permutations determine the map. The archive indexes the canonical codes
by their SHA-256 digests.

\newpage
\section{Geometric exclusion}

\begin{lemma}\label{lem:perturb}
An embedded surface in the stated class admits arbitrarily close
realisations of the same map in which every determinant of four distinct
plane vectors, including the plane at infinity, is nonzero.
\end{lemma}

\begin{proof}
Perturb the supporting planes independently and reconstruct each vertex
as the intersection of its incident triple. Transversality at the
vertices gives continuous dependence on the planes. Triangulate each
original face using its boundary vertices, and move the resulting
triangulation by moving its vertices. There are finitely many pairs of
simplices with disjoint vertex sets, and their images initially have
positive distance. They therefore remain disjoint under a sufficiently
small perturbation.

This condition ensures that the perturbed triangulation is embedded.
Indeed, if two distinct points of the abstract triangulation had the same
image, cancel their common barycentric weights and renormalise. The
remaining weights would exhibit an intersection between two simplices
with disjoint vertex sets. All triangles of a face remain in its perturbed
plane. The embedded triangulation therefore still represents a planar
polygon with the prescribed boundary. Finally, the required determinants
are finitely many nonzero polynomials in the plane coefficients. Their
simultaneous nonvanishing is dense among the small perturbations.
\end{proof}

Use homogeneous plane vectors \(h_i=(a_i,d_i)\), with equation
\(a_i\cdot x+d_i=0\), and \(h_\infty=(0,0,0,1)\).
The Boolean variables record the signs of the determinants of four
distinct vectors. The three-term Grassmann--Pl\"ucker relations impose
necessary conditions on these signs. This use of SAT to exclude oriented
matroids follows Schewe \cite{schewe}. The dual formulation in terms of
orders along intersections of face planes follows Moki and Julio
\cite{mokijulio}.

On the line \(L_{ij}\) of the planes \(i,j\), each remaining plane \(k\)
determines a point \(x_{ijk}\). With direction \(a_i\times a_j\),
\(x_{ij\ell}\) follows \(x_{ijk}\) exactly when
\[
 \operatorname{sgn}\det(h_i,h_j,h_\ell,h_k)\,
 \operatorname{sgn}\det(h_i,h_j,h_\ell,h_\infty)\,
 \operatorname{sgn}\det(h_i,h_j,h_k,h_\infty)=+1.
\]
These comparisons give the segment-membership tests used in the encoding.
The basic formula forbids intersections of nonconsecutive edges of a face
and intersections of the interior of an edge with a nonincident face.
Membership in a polygon is computed by the parity of boundary crossings,
with collinear boundary segments treated by the sides of their entering
and leaving edges. The strengthened formula also excludes an interval of
\(L_{ij}\) lying in both face interiors between allowed contacts. This
test matters when a pair of faces has several common edges. An overlap
along their line need not produce a transverse edge piercing.
Changing signs of the plane equations and reversing ambient orientation
allow redundant sign choices to be fixed without imposing any geometric
symmetry.

By Lemma~\ref{lem:perturb}, a realisation gives a sign assignment satisfying
these necessary conditions. A verified UNSAT certificate therefore
excludes a realisation. The basic formula suffices for \(39\) exclusions
in the eight-face face-neighbourly class. The other exclusions use the
strengthened formula. The verifier reconstructs the corresponding CNF,
checks its digest and runs DRAT-trim on the stored proof \cite{drat}.
The surviving sign assignments are resolved by exact geometric checks.

\newpage
\section{Counts and realisations}

\begin{table}[htbp]
\centering
\begin{tabular}{@{}rrlrrr@{}}
\toprule
\(F\) & \(g\) & Scope & Maps & UNSAT & Realised \\
\midrule
8 & 2 & Face-neighbourly & \(2\) & \(2\) & \(0\) \\
8 & 3 & Face-neighbourly & \(42\,122\) & \(42\,120\) & \(2\) \\
9 & 3 & Face-neighbourly & \(203\) & \(203\) & \(0\) \\
8 & 3 & Profile \(9^8\), unrestricted adjacency & \(103\) & \(100\) & \(3\) \\
\bottomrule
\end{tabular}
\caption{Isomorphism classes of maps and geometric outcomes. The last row
includes the \(51\) face-neighbourly maps of profile \(9^8\).}
\label{tab:counts}
\end{table}

The three genus-two profiles with eight faces are \(7^4 8^4\),
\(7^5 8^2 9\) and \(7^6 9^2\). They yield respectively one, zero and one
map. In the nine-face genus-three case the profile \(8^7 11^2\) is excluded
before enumeration. After reserving one edge per pair, only the two
11-gons have residual degree, equal to three each. The pair can carry at
most two additional edges.

The two surviving face-neighbourly maps match the coordinate realisations
in \cite{mizhaev2026,rostvigh}. Of the \(52\) maps added when face-neighbourliness
is dropped in profile \(9^8\), \(51\) are excluded and one is \(\Vone\).
In \(\Vone\), four pairs of faces have no common edge, twelve share one
edge, and twelve share two. For all three surviving maps, rational
coordinate data and cyclic face walks give exact checks of planarity,
simple boundaries and the prescribed intersections for all \(28\) face
pairs. This supplies the existence assertions of
Theorem~\ref{thm:classification}.

Interactive models of all three polyhedra are available on Geometheca.\footnote{%
Mizhaev's \(\Vone\) at \url{https://geometheca.org/polyhedra/mizhaev-octahedron-v1/}.\newline
Mizhaev's \(\Vtwo\) at \url{https://geometheca.org/polyhedra/mizhaev-octahedron-v2/}.\newline
The R\"ost--V\'igh polyhedron at \url{https://geometheca.org/polyhedra/rost-vigh-octahedron/}.}
The deposit uses the integer coordinates and cyclic face walks of the
Geometheca model for \(\Vone\). This is a realisation of the same map as
\(\Vone\) in \cite{mizhaev2020}, with different vertex coordinates.

\section*{Data and reproducibility}

The research deposit is available on Zenodo
(\href{https://zenodo.org/records/24122257}{DOI 10.5281/zenodo.24122257}).
It contains the profile lists, the maps, the SAT evidence, the DRAT proofs,
the coordinate data for the three realisations and the verification tools.
After extracting each map archive, run \verb|python tools/verify_research.py .|
from its root. The \verb|--enumerate| option also repeats map enumeration
using labelled triples and edge matchings and compares the canonical
signatures. This check addresses the coverage of the map list, while the
DRAT checks address the exclusion of individual maps.

\section*{Acknowledgements}

The author used OpenAI's ChatGPT to assist with translation, manuscript
formatting, writing Python code and carrying out some computations.
The author takes full responsibility for the content and correctness of this work.

\end{document}